\input amstex\input amsppt.sty
\magnification=\magstep1
\advance\vsize0cm\voffset=-1cm\advance\hsize1cm\hoffset0cm
\NoBlackBoxes
\define\R{{\Bbb R}}\define\Z{{\Bbb Z}}
\def\pr{\mathop{\fam0 pr}}

\def\forg{\mathop{\fam0 forg}}

\def\Sq{\mathop{\fam0 Sq}}
\def\diag{\mathop{\fam0 diag}}

\def\id{\mathop{\fam0 id}}

\def\rel#1{\allowbreak\mkern8mu{\fam0rel}\,\,#1}
\def\Int{\mathop{\fam0 Int}}

\def\Cl{\mathop{\fam0 Cl}}

\def\link{\mathop{\fam0 link}}
\def\im{\mathop{\fam0 im}}

\def\t{\widetilde}
\def\hd{1.1 }
\def\Hd{1.1}
\def\cor{1.2\ }
\def\Cor{1.2}
\def\wh{1.3\ }

\def\bn{1.4\ }

\def\co{2.1\ }
\def\Co{2.1}

\def\Wi{2.2}

\def\em{2.3\ }
\def\Em{2.3}
\def\di{2.4\ }
\def\Di{2.4}
\def\cs{2.5\ }
\def\Cs{2.5}
\def\tw{2.6\ }
\def\Tw{2.6}
\def\st{2.7\ }
\def\St{2.7}

\def\Pc{2.8}
\def\md{2.9\ }
\def\Md{2.9}
\def\tk{2.10\ }

\topmatter
\title Embeddings of $k$-connected $n$-manifolds into $\R^{2n-k-1}$ \endtitle
\author A. Skopenkov \endauthor
\address Department of Differential Geometry, Faculty of Mechanics and
Mathematics, Moscow State University, Moscow, Russia 119992, and the
Independent University of Moscow, B. Vlas\-yev\-skiy, 11, 119002, Moscow,
Russia. e-mail: skopenko\@mccme.ru \endaddress
\thanks
I would like to acknowledge D. Crowley, S. Melikhov, D. Tonkonog
and the anonymous referee for useful discussions and suggestions.
The author gratefully acknowledges the support from Deligne 2004 Balzan
prize in mathematics, the Russian Foundation for Basic Research Grants
07-01-00648a and 06-01-72551-NCNILa, President of Russian Federation Grant
MD-4729.2007.1
\endthanks

\abstract
We obtain geometric estimations for isotopy classes of embeddings of closed
$k$-connected $n$-manifolds into $\R^{2n-k-1}$ for $n\ge2k+6$ and $k\ge0$.
This is done in terms of an exact sequence involving the Whitney invariants
and an explicitly constructed action of $H_{k+1}(N;\Z_2)$
on the set of embeddings.
(For $k\ne1$ classification results were obtained by algebraic methods
without direct construction of embeddings or homology invariants.)

The proof involves a reduction to the classification of embeddings of a
punctured manifold and uses {\it the parametric connected sum} of embeddings.

{\bf Corollary.} {\it Suppose that $N$ is a closed almost parallelizable
$k$-connected $n$-manifold and $n\ge2k+6\ge8$.
Then the set of isotopy classes of embeddings $N\to\R^{2n-k-1}$
is in 1-1 correspondence with $H_{k+2}(N;\Z_2)$ for $n-k=4s+1$.}
\endabstract
\endtopmatter

\document
\head 1. Introduction and main results \endhead

{\bf Introduction.}

This paper
\footnote{An earlier version is published as [Sk10'].}
is on the classical Knotting Problem: {\it for an $n$-manifold $N$
and a number $m$ describe the set $E^m(N)$ of isotopy classes of embeddings
$N\to\R^m$}.
For recent surveys, see [RS99, Sk08, HCEC]; whenever possible we refer to these
surveys, not to original papers.

Denote CAT = DIFF (smooth) or PL (piecewise linear).
If the category is omitted, then a statement is correct (or a definition is
given) for both categories.

In this section assume that $N$ is a closed $k$-connected orientable $n$-manifold.

By $\Z_{(s)}$ we denote $\Z$ for $s$ even and $\Z_2$ for $s$ odd.

The Haefliger-Zeeman Unknotting Theorem [Sk08, Theorem 2.8.b] states that
{\it each two embeddings $N\to\R^m$ are isotopic for $n\ge2k+2$ and}
$$m\ge2n-k+1.$$
The classification of embeddings of $N$ into $\R^{2n-k}$ is as follows:
\footnote{This is a classical result of Haefliger-Hirsch (in the smooth
category) and Weber-Hudson-Bausum-Vrabec (in the PL category) [Sk08, Theorem
2.13].}
{\it for $n\ge2k+4$ and $k\ge0$ there is a 1--1 correspondence
(defined in Definition \wh below)}
$$W_{2n-k}: E^{2n-k}(N)\to H_{k+1}(N;\Z_{(n-k-1)}).$$
This paper is on the classification of
$$E^{2n-k-1}(N)$$
for $n\ge2k+6$.
Such a classification was known for $k=0$ [Ya84]
\footnote{I apologize for not presenting this reference in earlier versions of this
paper and I am grateful to S. Melikhov for bringing it to my attention in October, 2010.
Cf. [Ba75] for $N=\R P^n$ and [Hu69, BH70, Vr77, CS10, Sk10] for $n=4$.}
and for $k\ge2$ [BG71, Corollary 1.3].
\footnote{Namely, in [BG71, Corollary 1.3] it is stated that
{\it if $N$ is a closed $k$-connected orientable $n$-manifold embeddable
into $\R^m$, $m\ge2n-2k+1$
and $2m\ge3n+4$, then there is a 1--1 correspondence
$E^m(N)\to[N_0;V_{M,M+n-m+1}]$, where $M$ is large.}
We have $2n-k-1\ge 2n-2k+1$ if and only if $k\ge2$.
The exact sequences of Main Theorem \Hd.a below could apparently be obtained
using this result and homotopy classification of maps from an
$(n-k-1)$-polyhedron to an $(n-k-2)$-connected space. Cf. \S3.g.
}
Methods of [Ya84] possibly work for $k\ge1$; for an idea how to extend methods of [BG71] for $k=0,1$ see [To].
Those results were obtained using a reduction of the
classification of embeddings to an equivariant homotopy problem [Sk08, \S5].
They did not involve direct geometric constructions of embeddings or homology invariants.

In this paper we present such constructions (definitions of maps $b,W,W'$
in \S1, Remark \md in \S2 and a description of generators and relations of
$E^{2n-k-1}(S^{k+1}\times S^{n-k-1})$ in \S3).
Our proof is a generalization of the
Haefliger-Hirsch-Hudson-Vrabec argument for the proof of the bijectivity of
$W_{2n-k}$ (cited above).
Our constructions work for $k=1$ without any change.
Our result is explicit enough to yield Corollary \cor below.

\bigskip
{\bf Main results.}


An $n$-manifold $N$ is called {\it $p$-parallelizable} if each embedding
$S^p\to N$ extends to an embedding $S^p\times D^{n-p}\to N$.
\footnote{Note that 1-parallelizability is equivalent to orientability;
2-parallelizability is equivalent to 3-paralle\-li\-za\-bility and to the property
of being a spin manifold;
for each $p=4,5,6,7$ being $p$-parallelizable is equivalent to being a string
manifold.
A reader who is bothered by new terms can replace in this paper
the $p$-parallelizability by the almost parallelizability. }
If the coefficients of a homology group are omitted, then they are $\Z$.
For a group $G$, denote by $G*\Z_2$ the set of elements of order at most 2
in $G$.

\smallskip
{\bf Main Theorem \Hd.}
{\it (a) Let $N$ be a closed $k$-connected $n$-manifold.
Suppose that $k\ge1$, \ $n\ge2k+6$ and $N$ embeds
\footnote{The embeddability into $\R^{2n-k-1}$ is equivalent to
$\overline W_{n-k-1}(N)=0$, where $\overline W_{n-k-1}(N)$ is the normal
Stiefel-Whitney class [Sk08, \S2, Pr07, 11.3].}
into $\R^{2n-k-1}$.
For $n-k$ odd, assume that $N$ is $(k+2)$-parallelizable.
Then there is an exact sequence of sets
\footnote{The right-hand term can be represented by a
formula valid for both odd and even $n-k$:
\linebreak
$H_{k+2}(N)\otimes\Z_{(n-k)}\times H_{k+1}(N;\Z_{(n-k-1)})*\Z_2$.
The validity for $n-k$ even is obvious and for $n-k$ odd follows by the
Universal Coefficients Formula.}
with an action $b$:}
$$H_{k+1}(N;\Z_2)\overset b\to\to E^{2n-k-1}(N)
\cases
\overset{W\times W'}\to\to H_{k+2}(N)\times H_{k+1}(N;\Z_2)\to0
&n-k\text{ even,}\\
\overset W\to\to H_{k+2}(N;\Z_2)\to0 &n-k\text{ odd.}\endcases$$

\smallskip
{\bf Corollary \Cor.}
{\it Under the assumptions of Main Theorem \Hd.a

(a) for $n-k=4s+1$ there is a 1--1
correspondence
$$W: E^{2n-k-1}(N)\to H_{k+2}(N;\Z_2).$$
\quad
(b) the set $E^{2n-k-1}(N)$ is infinite if and only if $n-k$ is even
and $H_{k+2}(N)$ is infinite.}

\smallskip
Corollary \Cor.a is proved in \S2.
Corollary \Cor.b is straightforward.

For a closed connected $n$-manifold $N$ denote $N_0:=N-\Int B^n$, where
$B^n\subset N$ is a codimension 0 ball.
Consider the coefficient exact sequence
$$H_{k+2}(N)\overset2\to\to H_{k+2}(N)\overset{\rho_2}\to
\to H_{k+2}(N;\Z_2)\overset\beta\to\to H_{k+1}(N)\overset2\to\to H_{k+1}(N).$$
Here 2 is the multiplication by 2, $\rho_2$ is the reduction modulo 2 and
$\beta$ is the Bockstein homomorphism.

\smallskip
{\bf Main Theorem \Hd.}
{\it (b) Let $N$ be a closed connected orientable $n$-manifold.
\footnote{Each orientable $n$-manifold embeds into $\R^{2n-1}$
[Sk08, Theorem 2.4.a].}
If $n$ is odd, assume that $N$ is spin and the Hurewicz homomorphism
$\pi_2(N)\to H_2(N)$ is epimorphic.
For $n\ge6$ (and for $n=5$ in the PL category) there is an exact sequence of
sets with an action $b$
$$H_1(N;\Z_2)\overset b\to\to E^{2n-1}(N)\overset{W\times r}
\to\to\cases H_2(N)\times  E^{2n-1}(N_0)\to0 &n\text{ even}\\
H_2(N;\Z_2)\times E^{2n-1}(N_0)\overset a\to\to H_1(N)
&n\text{ odd}\endcases.$$
Here $r$ is the restriction-induced map and $a(x,f):=W_0'(f)-\beta(x)$.}

\smallskip
In Main Theorem \Hd.b the right-hand exactness implies that for $n$ odd $\im(W\times r)$ is in
1--1 correspondence with $\im\rho_2\times(2W_0')^{-1}(0)$. 


\smallskip
{\bf Open Problem.} {\it
(a) Find the preimages of $b$.
(See a discussion in \S3.e.)

(b) Describe the set $E^{2n-1}(N_0)$, cf. [Sa99, To], \S3, the Deleted Product Lemma.

(c) For $n\ge 2k+6$ there is a group structure on $ E^{2n-k-1}(N)$
(and for $n\ge4$ on $E^{2n-1}(N_0)$)
defined via the Haefliger-Wu $\alpha$-invariant [Sk08, \S5].
Are the maps from Main Theorem \hd homomorphisms?
If yes, solve the extension problem.}

\smallskip
{\bf Definitions \wh of the Whitney invariants $W$, $W_{2n-k}$, $W'$ and
$W_0'$.}
We present definitions for $n\ge2k+5$, $N$ orientable and
in the smooth category.
The definition in the PL category is analogous [Sk08, \S2.4].
\footnote{In Main Theorem \Hd.a,b, $k\ge1$, so $N$ is orientable.
For an equivalent definition see the Difference Lemma \di below or
[Sk08', \S1].}
Fix orientations on $N$ and on $\R^{2n-k-1}$.
Take embeddings $f,f_0:N\to\R^{2n-k-1}$.

The {\it self-intersection set} of a map $H:X\to Y$ is
$\Sigma(H):=\{x\in X\ |\ \#H^{-1}Hx>1\}.$

Take a general position homotopy $H:N\times I\to \R^{2n-k-1}\times I$ between
$f_0$ and $f$.
Since $n\ge2k+5$, by general position, $\Sigma(H)$ is a $(k+2)$-submanifold
(not necessarily compact).
The closure $\Cl\Sigma(H)$ is a closed $(k+2)$-submanifold.
For $n-k$ even it has a natural orientation.
\footnote{
{\it Definition of the orientation} is analogous to [Sk08, \S2.3, p. 263].
Take smooth triangulations $T$
and $T'$ of the domain and the range of $H$
such that $H$ is simplicial.
Then $\Cl\Sigma(H)$ is a subcomplex of $T$.
Take any oriented simplex $\sigma\subset\Cl\Sigma(H)$.
Let us show how to decide whether the orientation of $\sigma$ is right or to be
changed.
By general position there is a unique simplex $\tau$ of $T$ such that
$f\sigma=f\tau$.
The orientation on $\sigma$ induces an orientation on $f\sigma$ and then on
$\tau$.
The orientations on $\sigma$ and $\tau$ induce orientations on normal spaces in
$N\times I$ to these simplices.
These two orientations (in this order) together with the orientation on
$f\sigma$ induce an orientation on $\R^{2n-k-1}\times I$.
If this orientation agrees with the fixed orientation of $\R^{2n-k-1}\times I$,
then the orientation of $\sigma$ is right, otherwise it should be changed.
Since $n-k$ is even, these orientations agree for adjacent simplices
[Hu69, Lemma 11.4].
So they define an orientation of $\Cl\Sigma(H)$.
}
Define the Whitney invariant
$$W:E^{2n-k-1}(N)\to H_{k+2}(N\times I;\Z_{(n-k)})\cong H_{k+2}(N;\Z_{(n-k)})$$
by
$W(f)=W_{f_0}(f):=[\Cl\Sigma(H)].$
Analogously to [Sk08, \S2.4], this is well-defined.

The Whitney invariant
$W_{2n-k}: E^{2n-k}(N)\to H_{k+1}(N;\Z_{(n-k-1)})$
is defined analogously to the above.
The Whitney invariant $W': E^{2n-k-1}(N)\to H_{k+1}(N;\Z_{(n-k-1)})$
is defined as the composition of $W_{2n-k}$ and the map
$ E^{2n-k-1}(N)\to E^{2n-k}(N)$ induced by the inclusion
$\R^{2n-k-1}\to \R^{2n-k}$.

By [Vr89, Theorem 3.1] the map $W'$ equals (up to sign for $n-k$ odd)
the composition
$$ E\phantom{}^{2n-k-1}(N)\overset r\to\to E\phantom{}^{2n-k-1}(N_0)
\overset{W'_0}\to\to H_{k+1}(N;\Z_{(n-k-1)}).$$
Here $r$ is the restriction map and $W_0'$ is defined as follows.

The {\it singular set} of a smooth map $H:X\to Y$ between manifolds is
\linebreak
$S(H):=\{x\in X\ :\ d_xH\text{ is degenerate}\}$.

Take a general position homotopy $H:N\times I\to \R^{2n-k-1}\times I$ between
$f_0$ and $f$.
Since $n\ge2k+3$, by general position, $\Cl S(H)$ is a closed
$(k+1)$-submanifold.
For $n-k$ odd it has a natural orientation.
\footnote{{\it Definition of the orientation.}
Recall the notation from the previous footnote.
Take a $(k+1)$-simplex $\alpha\subset S(H)$.
By general position there are $(k+2)$-simplices
$\sigma,\tau\subset\Cl\Sigma(H)$ such that $f\sigma=f\tau$ and
$\sigma\cap\tau=\alpha$.
Define the `right' orientation of $\alpha$ to be the orientation induced by the
`right' orientation of $\sigma$.
This is well-defined because the `right' orientation of $\tau$ induces the same
orientation of $\alpha$.
(Indeed, since $n-k$ is odd, normal spaces of $\sigma$ and of $\tau$ in
$N\times I$ are even-dimensional, so the `right' orientations on $\sigma$ and
on $\tau$ induce the same orientation on $f\sigma$.)
}
Define $W'_0(f)=W'_{0,f_0}(f):=[S(H)]$.
It is well-known that $W'_0(f)$ is indeed independent of $H$ (for fixed $f$
and $f_0$).
\footnote{We use $W'_0r$, not $W'$ in the proof.
Although we do not need this, note that $W'_0r$ is a regular
homotopy invariant;
\quad
if $k=0$ and $n$ is even, then $W_0'$ factors through $H_{k+1}(N)$;
\quad
for $n-k$ odd, $2W_0'(f)$ equals the normal Euler class $e$ of $f$
because $e=AD^{-1}[f(\partial N_0)]=2W_0'(f)$ [Vr89, Addendum 2.2].}


\smallskip
For an embedding $f:N\to\R^{2n-k-1}\subset S^{2n-k-1}$ of a closed connected
$n$-manifold $N$ denote

$\bullet$ $C=C_f:=S^{2n-k-1}-f\Int N_0$,

$\bullet$ by $AD=AD_{f,i}:H_i(N)\to H_{i+n-k-2}(C)$ the composition of
Alexander and Poincar\'e isomorphisms,

$\bullet$ by $h=h_{f,i}:\pi_i(C)\to H_i(C)$ the Hurewicz homomorphism.

\smallskip
{\bf Definition \bn of the action $b=b_N$.}
\footnote{A reader who is not interested in explicit constructions can omit
this definition and set $b(x)f:=\psi_fb'(x)$, where $b'$ and $\psi_f$ are
defined in \S2.}
We give a definition for $n\ge2k+6$ and $H_k(N)=0$.
Take an embedding $f:N\to \R^{2n-k-1}$ and $x\in H_{k+1}(N;\Z_2)$.
Since $H_k(N)=0$, there is $\overline x\in H_{k+1}(N)$ such that
$\rho_2\overline x=x$.
By general position and Alexander duality, $C$ is $(n-2)$-connected.
Hence $h_{n-1}$ is an isomorphism.
Consider the composition
$S^n\overset{\Sigma^{n-3}\eta}\to\to S^{n-1}
\overset{h_{n-1}^{-1}AD\overline x}\to\to C,$
where $\eta:S^3\to S^2$ is the Hopf map.
\footnote{
Since $N$ is $k$-connected and $n\ge2k+3$, we can represent $\overline x$ by an
embedding $x':S^{k+1}\to N$.
If the restriction to $x'(S^{k+1})$ of the normal bundle
$\nu_f:\partial C\to N$ is trivial, then {\it spheroid
$h_{n-1}^{-1}AD\overline x$ can be constructed directly as follows}, cf.
[Sk08', end of \S1].
Identify $X:=\nu_f^{-1}x'(S^{k+1})$ with $S^{k+1}\times S^{n-k-2}$.
Let us show how to make an embedded surgery of $S^{k+1}\times*\subset X$ to
obtain an $(n-1)$-sphere $S^{n-1}\cong\Sigma\subset C$ whose inclusion into $C$
represents $h_{n-1}^{-1}AD\overline x$.
\newline
Take a vector field on $S^{k+1}\times*$ normal to $X$ in $\R^{2n-k-1}$.
Extend $S^{k+1}\times*$ along this vector field to a smooth map
$\widetilde x:D^{k+2}\to S^{2n-k-1}$.
Since $2n-k-1>2k+4$ and $n+k+2<2n-k-1$, by general position we may assume that
$\widetilde x$ is a smooth embedding and $\widetilde x(\Int D^{k+2})$ misses
$f(N)\cup X$.
Denote $l:=2n-2k-3$.
Since $n-k-1>k+1$, we have $\pi_{k+1}(V_{l,n-k-2})=0$.
Hence the standard framing of $S^{k+1}\times*$ in $X$
extends to an $l$-framing on $\widetilde x(D^{k+2})$ in $\R^{2n-k-1}$.
Thus $\widetilde x$ extends to an embedding
$\widehat x:D^{k+2}\times D^l\to C$ such that
$\widehat x(\partial D^{k+2}\times D^l)\subset X.$
Let
$$\Sigma:\ =\ (X-\widehat x(\partial D^{k+2}\times\Int D^l))
\bigcup\limits_{\widehat x(\partial D^{k+2}\times\partial D^l)}
\widehat x(D^{k+2}\times\partial D^l)\ \cong\ S^{n-1}.$$
}
The connected sum in $C$ of this composition with $f|_{B^n}$ is homotopic
(relative to the boundary) to
an embedding $x'':B^n\to C$ by Theorem \em below.
Define $b(x)f$ to be $f$ on $N_0$ and $x''$ on $B^n$.

This is well-defined (i.e. is independent of the choices of $\overline x$ and
of $x''$) for $n\ge2k+6$ and is an action by the equivalent definition given
in the proof of the Construction Lemma \cs below.

\head 2. Proof of Main Theorem \Hd \endhead

{\bf Main tools.}

The proof is based on the construction and application of the following
commutative diagram:
$$\CD H_{k+1}(N;\Z_2) @>> b> @. @. \\
@VVb' V  @VV b V @. @.\\
\pi_n(C) @>>\psi_f> r^{-1}r(f) \subset E^{2n-k-1}(N) @>>r>
E^{2n-k-1}(N_0)\\
@VV AD^{-1}\circ h V
@VV W V @VV W_0' V\\
H_{k+2}(N) @>>\rho_{(n-k)}>
H_{k+2}(N;\Z_{(n-k)}) @>>\beta >H_{k+1}(N;\Z_{(n-k-1)})\\
@VV V@.@.\\
0@.@.@.@.
\endCD.$$
Here

$\bullet$ $N$ is a closed homologically $k$-connected orientable $n$-manifold,
$f:N\to\R^{2n-k-1}$ is an embedding,  $n\ge2k+6$ and $N_0:=N-\Int B^n$, where
$B^n\subset N$ is a codimension 0 ball,

$\bullet$  $C$, $AD$, $h$ are defined at the end of \S1,

$\bullet$ $W$ and $W_0'$ are defined above in Definitions \wh of the Whitney
invariants,

$\bullet$ $r$ is the restriction-induced map,

$\bullet$ $\beta$ is the Bockstein homomorphism defined only for $n-k$ odd,

$\bullet$ $\psi_f$ is defined below in the Construction Lemma \Cs,

$\bullet$ $\rho_{(n-k)}$ is the identity for $n-k$ even and
the reduction modulo 2 for $n-k$ odd,

$\bullet$ $b$ is defined in \S1 (and can be alternatively defined as
$b:=\psi_f b'$),

$\bullet$ $b'$ is the composition
$$H_{k+1}(N;\Z_2)\overset\cong\to\to H_{n-1}(C;\Z_2)\overset\cong\to
\to H_{n-1}(C)\otimes\Z_2\overset\cong\to\to
\pi_{n-1}(C)\otimes\pi_n(S^{n-1})\to\pi_n(C)$$
of the Alexander duality, the coefficient isomorphism, tensor product of
the Hurewicz and the Pontryagin isomorphisms, and the composition map.
\footnote{The composition $\pi_{n-1}(C)\times\pi_n(S^{n-1})\to\pi_n(C)$ is
clearly linear in $\pi_n(S^{n-1})$; for $n\ge4$ the composition is linear in
$\pi_{n-1}(C)$ by [Po85, Lecture 4, Corollary in p. 167].
So the latter composition map is indeed well-defined.}

The proof of Main Theorem \hd in the next subsection shows how to apply this
diagram.
That proof uses statements of lemmas below not their proofs.

\smallskip
{\bf Complement Lemma \Co.}
{\it Let $N$ be a closed homologically $k$-connected orientable $n$-manifold,
$n\ge4$ and $f:N\to\R^{2n-k-1}$ an embedding.
Then the left column of the above diagram is exact. }

\smallskip
{\it Proof.}
By general position and Alexander duality, $C$ is $(n-2)$-connected.
Since $n\ge4$, by [Wh50] there is an exact sequence
forming the first line of the following diagram:
$$\minCDarrowwidth{0pt}\CD
H_{n-1}(C;\Z_2)@>>b'\circ AD^{-1}> \pi_n(C)@>>h> H_n(C)@>> >0\\
@AA {AD\cong} A  @.  @AA {AD\cong} A @.\\
H_{k+1}(N;\Z_2)  @. @. H_{k+2}(N) @.\endCD.$$
Now the lemma follows by  Alexander duality.
\qed


\smallskip
{\bf The Whitney Invariant Lemma \Wi.}
{\it Let $N$ be a closed $k$-connected orientable $n$-manifold
embeddable into $\R^{2n-k-1}$.

($W_0'$) The map $W_0'$ is a 1--1 correspondence for $k\ge1$;

the map $W_0'$ is surjective for $k=0$ and $n$ even;

$\im W_0'\supset H_1(N)*\Z_2$ for $k=0$ and $n$ odd.

(r) If $n\ge2k+6$, then $r$ is surjective for $n-k$ even and
$\im r=\ker(2W_0')$ for $n-k$ odd.

($\beta$)
\footnote{This is analogous to the
well-known
relations $w_{2j+1}=\Sq^1w_{2j}$ and $W_{2j+1}=\beta w_{2j}$ for the
Stiefel-Whitney classes [Pr07, 11.3]. }
For $n-k$ odd, $W_0'r=\beta W$. }

\smallskip
Part ($W_0'$) for $k=0$ follows by [Ya83, Main Theorem (i) and (iii)].
Part ($W_0'$) for $k\ge1$ and part ($r$) are proved in the PL category in
[Vr89, Theorem 2.1, Theorem 2.4 and Corollary 3.2] and in the smooth category
in [Ri70].
For the reader's convenience, the proofs of ($W_0'$) and ($r$) are sketched
below.

\smallskip
{\it Sketch of the proof of ($W_0'$).}
Let $Y$ be the set of regular homotopy classes of immersions
$N_0\to\R^{2n-k-1}$.
Since $N$ is $k$-connected, $N_0$ collapses to an $(n-k-1)$-polyhedron.
So by general position the forgetful map $ E^{2n-k-1}(N_0)\to Y$ is
surjective and, for $k\ge1$, injective (see details e.g. in [Vr89, proof of
Theorem 2.1] on p. 167).
The map $W_0'$ is a composition of the forgetful map and a certain map
$Y\to H_{k+1}(N;\Z_{(n-k-1)})$ that is a 1--1 correspondence for $k\ge1$
by the Smale-Hirsch
(in the smooth category)
or the Haefliger-Poenaru
(in the PL category)
classification of
immersions.

Now assume that $k=0$.
Then $\im W_0'\supset\im W'$.
By [Ya83, Main Theorem (i) and (iii)] $W'$ is surjective for $n$ even
and $\im W'=\im\beta=H_1(N)*\Z_2$ for $n$ odd.
This implies the required result on $\im W_0'$.
\qed
\footnote{For $k=0$ the map $W_0'$ is not injective.
For $k=0$ and $n$ odd, $\im W_0'$ can be larger than $H_1(N)*\Z_2$ (because
not all embeddings $D^1\times S^{n-1}\to\R^{2n-1}$ extend to embeddings
$S^1\times S^{2n-1}\to\R^{2n-1}$).
The result [Vr89, Corollary 3.1] holds for $k=0$ and orientable $N$
[Ya83, Main Theorem (2.i,iii)].}

\smallskip
{\bf Theorem \Em.} [RS99, Theorem 3.2]
{\it Let $N$ and $M$ be $n$- and $m$-manifolds with boundary.
Assume that
$2m\ge 3n+4$.

(a) If $N$ is $(2n-m)$-connected and $M$ is $(2n-m+1)$-connected, then any
proper map $N\to M$ whose restriction to the boundary $\partial N$ is an
embedding is homotopic (relative to the boundary $\partial N$) to an embedding.

(b) If $N$ is $(2n-m+1)$-connected and $M$ is $(2n-m+2)$-connected, then any
proper homotopy $N\times I\to M\times I$ fixed on the boundary $\partial N$ is
homotopic (relative to $\partial(N\times I)$) to an isotopy.}

\smallskip
{\it Sketch of the proof of ($r$).}
Let $f:N_0\to\R^{2n-k-1}$ be an embedding.
If $n-k$ is even, then the {\it homology} class of $f(\partial N_0)$ in
$H_{n-1}(C)$ is trivial [Vr89, proof of Theorem 2.4 and Addendum 2.2].
By general position and Alexander duality $C$ is $(n-2)$-connected.
Hence $h_{n-1}$ is an isomorphism.
Therefore the {\it homotopy} class of $f(\partial N_0)$ in $\pi_{n-1}(C)$
is trivial.
Then by Theorem \Em.a $f$ extends to an embedding $N\to\R^{2n-k-1}$.
Thus $r$ is surjective.

If $n-k$ is odd, then the homology class of $f(\partial N_0)$ equals
$2AD(W'_0(f))$ [Vr89, proof of Theorem 2.4 and Addendum 2.2].
Thus $\im r=\ker(2W_0')$ analogously to the case when $n-k$ is even.
\qed
\footnote{For $n-k$ odd the inclusion $\im r\subset\ker(2W_0')$ of part (r)
also follows by part ($\beta$) or by an analogue of the Boechat-Haefliger Lemma
[Sk08', \S2].}

\smallskip
{\it Proof of ($\beta$).}
Take a general position homotopy $H$ between $f_0$ and $f$.
Recall that in the Definition \wh of the Whitney invariants (including
footnotes) we defined integer $(k+2)$- and $(k+1)$-chains $[\Cl\Sigma(H)]$ and
$[\Cl S(H)]$ in $N\times I$
(simplicial chains in a certain smooth triangulation).
The assumption that $n-k$ is even was only used to show that
$\partial[\Cl\Sigma(H)]=0$; the assumption that $n-k$ is odd was used to define
$[\Cl S(H)]$.

Take two $(k+2)$-simplices $\sigma,\tau\subset \Cl\Sigma(H)$ intersecting by a
$(k+1)$-simplex $\alpha$.
Clearly, for $\Int\alpha\subset\Sigma(H)$ (the `right' orientations of
$\sigma$ and $\tau$ agree and) $\alpha$ appears in $\partial\sigma$ and in
$\partial\tau$ with the opposite signs.
Since $n-k$ is odd, for $\alpha\subset\Cl S(H)$ (the `right' orientations of
$\sigma$ and $\tau$ disagree and) $\alpha$ appears in $\partial\sigma$ and in
$\partial\tau$ with the same sign.
This and $\Cl\Sigma(H)=\Sigma(H)\cup S(H)$ imply that
$\partial[\Cl\Sigma(H)]=2[S(H)]$  for $n-k$ odd.

Then $\beta W(f)=[S(H)]=W_0'r(f)$.
Here the first equality holds by definition of $\beta$, and the second equality
holds by definition of $W_0'$.
\qed

\smallskip
{\bf Difference Lemma \Di.}
{\it Let $N$ be a closed connected orientable $n$-manifold and
$f,f':N\to\R^{2n-k-1}$ embeddings coinciding on $N_0$.
Then
\footnote{In this formula $\overline B^n$ is $B^n$ with reversed
orientation; we have $C_f=C_{f'}$.
The element $d(f',f)$ is an invariant of an isotopy (of $f$ and $f'$) relative
to $N_0$. }
}
$$W(f)-W(f')=\rho_{(n-k)}d(f',f),\quad\text{where}
\quad d(f',f):=AD_f^{-1}h_{f,n}[f'|_{B^n}\cup f|_{\overline B^n}]
\in H_{k+2}(N).$$

{\it Proof.} Take a map $F:B^{n+1}\to\R^{2n-k-1}$ in general position with
$f(N_0)$ and such that $F|_{\partial B^{n+1}}=f'|_{B^n}\cup f|_{\overline B^n}$.
By Alexander duality, $d(f,f')$ is the homology class carried by
$f^{-1}F(\Int B^{n+1})$.
There is a general position homotopy $H$ between $f$ and $f'$ such that
$\pr_N\Cl\Sigma(H)=f^{-1}F(\Int B^{n+1})$.
For $n-k$ even observe that in this formula the signs of corrresponding
simplices (in a certain smooth or PL triangulation of $N$) are the same.
So the lemma follows.
\qed

\smallskip
{\bf Construction Lemma \Cs.}
{\it Let $N$ be a closed homologically $k$-connected orientable $n$-manifold,
$f:N\to\R^{2n-k-1}$ an embedding and $n\ge2k+6$.
Then there is a map
$\psi=\psi_f:\pi_n(C)\to r^{-1}r(f)$ such that

(a) $d(\psi(y),f)=AD^{-1}h_n(y)$ for each $y\in\pi_n(C)$.

(b) $W(\psi(y))-W(f)=\rho_{(n-k)}AD^{-1}h_n(y)$ for each $y\in\pi_n(C)$.

(c) If $f=f'$ on $N_0$, then $f'$ is isotopic to
$\psi[f'|_{B^n}\cup f|_{\overline B^n}]$ relative to $N_0$.

(d) $\psi$ is surjective.

(e) $\psi$ defines an action.}

\smallskip
{\it Proof.}
Construction of $\psi$ is analogous to [Sk08', proof of the surjectivity of $W$
in \S5].
Take $x\in\pi_n(C)$ represented by a map $x':S^n\to C$.
The connected sum $x'\#f|_{B^n}$ in $C$ of $x'$ with $f|_{B^n}$ is homotopic
$\rel\partial B^n$ to a proper embedding $x'':B^n\to C$ coinciding with $f$ on
$\partial B^n$,
and $x''$ is uniquely defined by $x$ up to isotopy $\rel\partial B^n$.
(This follows by Theorem \em because
by general position and Alexander duality $C$ is $(n-2)$-connected,
$2n-(2n-k-1)+2\le n-2$ and $2(2n-k-1)\ge3n+4$.)
Define $\psi(x)$ to be $f$ on $N_0$ and $x''$ on $B^n$.

Since $f=\psi(x)$ on $N_0$, we have $\psi(x)\in r^{-1}r(f)$.

Part (a) holds because $y=[\psi(y)|_{B^n}\cup f|_{\overline B^n}]$.

Part (a) and the Difference Lemma \di imply (b).

Part (c) follows analogously to the uniqueness of $x''$ in the
construction of $\psi$.

If $r(f_1)=r(f)$ for an embedding $f_1:N\to\R^{2n-k-1}$, then $f_1$ is isotopic
to an embedding $f'$ such that $f=f'$ on $N_0$.
Then by (c) $\psi[f'|_{B^n}\cup f|_{\overline B^n}]$ is isotopic $\rel N_0$ to $f'$ and hence to $f_1$.
This implies (d).

Let us prove
part (e).
Take $x,y,x+y\in\pi_n(C)$ represented by maps $x',y',(x+y)':S^n\to C$.
We have that $x'\#(y'\#f|_{B^n})$ is homotopic $\rel\partial B^n$ to
$(x'+y')\#f|_{B^n}$.
Hence $\psi_f(x+y)$ is isotopic $\rel N_0$ to $\psi_{\psi_f(y)}(x)$ analogously
to the uniqueness of $x''$ in the construction of $\psi$.
\qed

\smallskip
For $x\in H_{k+1}(N;\Z_2)$ define $b(x)f:=\psi_fb'(x)$.
(Recall that $b'$ is defined in the Complement Lemma \Co;
this is clearly equivalent to the definition given in \S1.)

\bigskip
{\bf Proof of Main Theorem \Hd.}

\smallskip
{\it Proof of Main Theorem \Hd.a for $n-k$ even.}
The map $W'=W_0'r$ is surjective by the Whitney Invariant Lemma \Wi.r,$W_0'$.
Since $\rho_{(n-k)}=\id$ and $h_n$ is epimorphic, by the Construction Lemma
\Cs.b,d, $W\times W'$ is surjective.

Let us prove the exactness at $ E^{2n-k-1}(N)$.

By the Complement Lemma \Co, $h_nb'=0$.
Hence by the Difference Lemma \Di,  $W(f)=W(\psi_fb'(x))$ for each $x$.
Since $r$ is a factor of $W'$ and $r(f)=r(\psi_fb'(x))$, we have
$W'(f)=W'(\psi_f b'(x))$ for each $x$.

Suppose that $W(f)=W(g)$ and $W'(f)=W'(g)$.
Then $r(f)=r(g)$ by the Whitney Invariant Lemma \Wi.$W_0'$ because $W'=W_0'r$.
Thus $g=\psi_f(y)$ for some $y\in\pi_n(C)$ by
the Construction Lemma \Cs.d.
By the Construction Lemma \Cs.b, $h_n(y)=0$.
Hence by the Complement Lemma \Co, $y=b'(x)$ for some $x\in H_{k+1}(N;\Z_2)$.
So $g=\psi_fb'(x)=b(x)f$.
\qed

\smallskip
{\it Proof of Main Theorem \Hd.b for $n$ even.}
By the Whitney Invariant Lemma \Wi.r and the Construction Lemma \Cs.b,d the
map $r\times W$ is surjective.

Clearly, $r(f)=r(\psi_fb'(x))$.
Analogously to the previous proof,

$\bullet$ $W(f)=W(\psi_f b'(x))$,

$\bullet$ if $W(f)=W(g)$ and $r(f)=r(g)$, then $g=\psi_f b'(x)$ for some
$x\in H_1(N;\Z_2)$.
\qed

\smallskip
Now we turn to the case when $n-k$ is odd.
The proof of the following result is postponed.

\smallskip
{\bf Twisting Lemma \Tw.}
{\it Suppose that $n-k$ is odd, $N$ is a closed connected
$(k+2)$-parallelizable $n$-manifold and $n\ge2k+6$.
Assume that the Hurewicz homomorphism $\pi_{k+2}(N)\to H_{k+2}(N)$ is
epimorphic (for $k\ge1$ this follows from the $k$-connectedness).
Then for each $x\in H_{k+2}(N)$ every embedding $f:N\to\R^{2n-k-1}$ is isotopic
to an embedding $f':N\to\R^{2n-k-1}$ such that $f=f'$ on $N_0$ and
$d(f',f)=2x\in H_{k+2}(N)$.}

\smallskip
{\it Proof of Main Theorem \Hd.a for $n-k$ odd.}
By the Whitney Invariant Lemma \Wi.r,$W_0'$,$\beta$ we have
$\im(\beta W)=\im(W_0'r)=H_{k+1}(N)*\Z_2=\im\beta$.
Since $h_n$ is surjective, by the Construction Lemma \Cs.b,d,
$W(r^{-1}r(f))=W(f)+\im\rho_2=W(f)+\ker\beta$.
Thus $W$ is surjective.


Analogously to the case of $n-k$ even $W(f)=W(\psi_fb'(x))$ for each $x$.
If $W(f)=W(g)$, then $W_0'r(f)=W_0'r(g)$ by the Whitney Invariant Lemma
\Wi.$\beta$.
Hence by the Whitney Invariant Lemma \Wi.$W_0'$ we have $r(f)=r(g)$,
i.e. $g$ is isotopic to an embedding $g_1$ such that $g_1=f$ on $N_0$.
Since $W(g_1)=W(g)=W(f)$, by the Difference Lemma \Di, $d(g_1,f)$ is even.
Hence by the Twisting Lemma \Tw, $g_1$ is isotopic to an embedding $g_2$
such that $g_2=f$ on $N_0$ and $d(g_2,f)=0$.
By
the Construction Lemma \Cs.d, there is $y\in\pi_n(C)$ such that $g_2$ is
isotopic to $\psi_f(y)$ relative to $N_0$.
By the Construction Lemma \Cs.a, $AD^{-1}h_n(y)=d(\psi_f(y),f)=d(g_2,f)=0$.
Hence by the Complement Lemma \Co, $y=b'z$ for some $z\in H_{k+1}(N;\Z_2)$.
Thus $g$ is isotopic to $\psi_fb'(z)=b(z)f$.
\qed

\smallskip
{\it Proof of Main Theorem \Hd.b for $n$ odd.}
If $W(f)=W(g)$ and $r(f)=r(g)$, then analogously to the proof of (a)
for $n-k$ odd, $g=b(z)f$ for some $z\in H_1(N;\Z_2)$.

By the Whitney Invariant Lemma \Wi.$\beta$ we have $W_0'r=\beta W$, so
$\im(W\times r)\subset\ker a$.

Let us prove that $\im(W\times r)\supset\ker a$.
Take $x\in H_2(N;\Z_2)$ and $f:N_0\to\R^{2n-1}$ such that $\beta(x)=W_0'(f)$.
Then $2W_0'(f)=2\beta(x)=0$.
Hence by the Whitney Invariant Lemma \Wi.r, $f$ extends to an embedding
$f_1:N\to\R^{2n-1}$.
By the Whitney Invariant Lemma \Wi.$\beta$ we have
$\beta W(f_1)=W_0'(f)=\beta(x)$.
Hence $W(f_1)-x=\rho_2y'$ for some $y'\in H_2(N)$.
Since $h_{f_1,n}$ is surjective, there is $y\in\pi_n(C_{f_1})$ such that
$AD_{f_1}h_{f_1,n}(y)=y'$.
Then by the Construction Lemma \Cs.b,
$$\psi_{f_1}(-y)=f_1=f\quad\text{on}\quad N_0\quad\text{and}
\quad W(\psi_{f_1}(-y))=W(f_1)-\rho_2y'=x.\quad\qed$$

\smallskip
{\bf Parametric connected sum of embeddings.}

In this subsection we recall, with only minor modifications, some results of
[Sk07], cf. [PCS].

Denote $D^k_\pm:=\{(x_0,x_1,\dots,x_k)\in S^k\ |\ \pm x_0\ge0\}$.
Identify $D^p$ with $D^p_+$ and $S^p$ with
\quad $D^p_+\bigcup\limits_{\partial D^p_+=\partial D^p_-}D^p_-.$

For $m\ge n+2$ denote by $t^m_{p,n-p}$ the CAT {\it standard embedding} that is the
composition
$$S^p\times S^{n-p}\to\R^{p+1}\times\R^{n-p+1}\to\R^m\to S^m$$
of CAT standard embeddings.

Take an embedding $s:S^p\times D^{n-p}\to N$.
A map $f:N\to S^m$ is called {\it $s$-standardized} if

$\bullet$ $f(N-\im s)\subset\Int D^m_+$ and

$\bullet$ $f\circ s:S^p\times D^{n-p}\to D^m_-$ is the
restriction of the standard embedding.

A map $F:N\times I\to S^m\times I$ such that $F|_{N\times j}:N\times j\to S^m\times j$
is $s$-standardized (for $j=0,1$) is called {\it $s$-standardized} if

$\bullet$ $F((N-\im s)\times I)\subset\Int D^m_+\times I$ and

$\bullet$ $F\circ(s\times\id I):S^p\times D^{n-p}\times\{t\}
\to D^m_-\times\{t\}$
is the restriction of the standard embedding for each $t\in I$.

\smallskip
{\bf Standardization Lemma \St.} [Sk07, Standardization Lemma]
{\it Let $N$ be an $n$-manifold $N$ and  $s:S^p\times D^{n-p}\to N$
an embedding.
For $m\ge n+p+3$,

$\bullet$ each embedding $N\to S^m$ is isotopic to an $s$-standardized
embedding, and

$\bullet$ each concordance between $s$-standardized embeddings is isotopic
relative to the ends to an $s$-standardized concordance.}

\smallskip
Denote $T^{p,n-p}:=S^p\times S^{n-p}.$
Let $i:S^p\times D^{n-p}\to T^{p,n-p}$ be the standard inclusion.
Recall that $N_0=N-\Int B^n$, where $B^n\subset N$ is a codimension 0 ball.

\smallskip
{\bf Summation Lemma \Pc.}
{\it Assume that  $m\ge n+p+3$, $N$ is a closed connected
$n$-manifold and $f:N\to S^m$, $g:T^{p,n-p}\to S^m$ embeddings.

(a) By the Standardization Lemma \st we can make concordances
and assume that $f$ and $g$ are $s$-standardized and $i$-standardized,
respectively.
Then an embedding
$$f\#_sg:N\to S^m\quad\text{is well-defined by}\quad
(f\#_sg)(a)=\cases f(a) &a\not\in\im s\\
R_mg(x,R_{n-p}y)           &a=s(x,y)\endcases,$$
where $R_k$ is the symmetry of $S^k$ with respect to the hyperplane
$x_1=x_2=0$.

(b) If $g=t^m_{p,n-p}$ on $(T^{p,n-p})_0$, then $f\#_sg=f$ on $N_0$.

(c) If $m=2n-p+1$ and $g=t^m_{p,n-p}$ on $(T^{p,n-p})_0$, then
$d(f\#_sg,f)=d(g,t^m_{p,n-p})[s|_{S^p\times 0}]\in H_p(N)$, where $p<n/2$ and
$d(g,t^m_{p,n-p})\in H_p(T^{p,n-p})$ is considered as an integer.}

\smallskip
{\it Proof.}
The argument for (a) is easy and similar to [Sk06, Sk07].
In order to prove that $f\#_sg$ is well-defined we need to show that
the concordance class of $f\#_sg$ depends only on concordance classes of
$f$ and $g$ but not on the chosen standardizations of $f$ and $g$.
Take concordances
$$F:N\times I\to S^m\times I\quad\text{and}
\quad G:T^{p,n-p}\times I\to S^m\times I$$
between different standardizations of $f$ and of $g$.
By the `concordance' part of the Standardization Lemma \st we can take
concordances relative to the ends and assume that $F$ and $G$ are
$s$-standardized and $i$-standardized, respectively.
Define a concordance
$$F\#_sG:N\times I\to S^m\times I\quad\text{by}
\quad(F\#_sG)(a,t)=\cases F(a,t)& a\not\in\im s\\
R_mG(x,R_{n-p}y,t)      &a=s(x,y)\endcases.$$
If $F$ is a concordance from $f_0$ to $f_1$ and $G$ is a concordance from $g_0$
to $g_1$, then $F\#_sG$ is a concordance from $f_0\#_sg_0$ to $f_1\#_sg_1$.

Parts (b) and (c) are clear.
\qed



\bigskip
{\bf Applications of the parametric connected summation.}

\smallskip
{\it Proof of the Twisting Lemma \Tw.}
Denote $t=t_{k+2,n-k-2}^{2n-k-1}$ and $t_1:=\psi_t(h_{t,n}^{-1}AD_t(2))$
for the generator $1\in H_{k+2}(T^{k+2,n-k-2})\cong\Z$
(the map $h_{t,n}$ is an isomorphism).
Since $n-k$ is odd, by the Construction Lemma \Cs.b,
$$W(t_1)-W(t)=\rho_{(1)}(2)=0\in H_{k+2}(T^{k+2,n-k-2};\Z_2)\cong\Z_2.$$
Hence $t_1$ is isotopic to $t$ by the bijectivity of $W_{2n-k}$
[Sk08, Theorem 2.13].

Since the Hurewicz homomorphism $\pi_{k+2}(N)\to H_{k+2}(N)$ is epimorphic
and $n\ge2k+5$, by general position there is an embedding $S^{k+2}\to N$
realizing $x$.
Since $N$ is $(k+2)$-parallelizable, this embedding extends
to an embedding $\overline x:S^{k+2}\times D^{n-k-2}\to N$.

Since $2n-k-1\ge n+k+2+3$, embeddings $f=f\#_{\overline x}t$ and
$f':=f\#_{\overline x}t_1$ are well-defined by the Summation Lemma \Pc.a
and are isotopic.
We have  $d(f\#_{\overline x}t_1,f)=d(t_1,t)x=2x$ by the Construction Lemma
\Cs.a and the Summation Lemma \Pc.c.
\qed

\smallskip
{\it Proof of Corollary \Cor.a.}
By Main Theorem \Hd.a it remains to prove that $b_N=0$ for $n-k=4s+1$.

Since $N$ is $k$-connected, the composition
$\pi_{k+1}(N)\to H_{k+1}(N)\overset{\rho_2}\to\to H_{k+1}(N;\Z_2)$
of the Hurewicz isomorphism and the reduction modulo 2 is an epimorphism.
Hence for each $x\in H_{k+1}(N;\Z_2)$ there is an embedding $S^{k+1}\to N$
realizing $x$ (because $n\ge 2k+3$).
Since $N$ is $(k+1)$-parallelizable, this embedding extends
to an embedding $\overline x:S^{k+1}\times D^{n-k-1}\to N$.

Denote $$T:=S^{k+1}\times S^{n-k-1},\quad t:=t^{2n-k-1}_{k+1,n-k-1}\quad\text{and}
\quad  \gamma:=b_T(1)t,$$
where $1\in H_{k+1}(T;\Z_2)$ is the generator.
Since $n-k\equiv1\mod4$, by [Sk08, Theorem 3.9 and tables] $\gamma$ is
isotopic to $t$.
Therefore
$$b_N=0\quad\text{because}
\quad b_N(x)f=f\#_{\overline x}\gamma=
f\#_{\overline x}t=f.$$
Here the parametric connected sums are well-defined because
$2n-k-1\ge n+k+1+3$.
In order to prove the first equality we assume in the construction of $\gamma$
that $S^n, S^{n-1}\subset \R^{2n-k-1}_+$.
Then we may assume that $\gamma$ is standardized.
So $f\#_{\overline x}\gamma$ is obtained from $f$ by linked connected
summation along $\overline x$ with a composition
$S^n\to S^{n-1}\times D^{n-k}\to S^{2n-k-1}-f(N)$ of
two embeddings, the one representing $\Sigma^{n-3}\eta$ and the other
representing $AD(x)$.
Hence $b_N(x)f=f\#_{\overline x}\gamma$ by definition of $b$.
\qed

\smallskip
{\bf Remark \Md.}
{\it Let $N$ be a closed $k$-connected $(k+2)$-parallelizable $n$-manifold,
$n\ge2k+6$ and $k\ge1$.
Then every embedding $N\to\R^{2n-k-1}$ can be obtained from every other
embedding by parametric connected summations with embeddings
($\gamma$ is defined in the above proof; $\tau$ and $\varkappa$ are defined
below):

$\bullet$ $\gamma$, $\varkappa$ and $\tau$, provided $n-k$ is even;

$\bullet$ $\gamma$ and $\varkappa$, provided $n-k$ is odd and $H_{k+1}(N)$ has
no 2-torsion.
\footnote{It would be interesting to drop the latter condition; for this, one
needs an explicit construction of embeddings whose Whitney invariants are in
$\im\beta\subset H_{k+2}(N;\Z_2)$.
For this, one needs an explicit construction of immersions $S^n\to\R^{2n}$.}
}

\smallskip
{\bf Definition \tk of embeddings $\tau$ and $\varkappa$.}
Define the Hudson Torus $\varkappa:T^{k+2,n-k-2}\to \R^{2n-k-1}$ as in
[Sk08, \S2.2] or set $\varkappa:=\psi_t(h_{t,n}^{-1}AD_t(1))$ for the
generator $1\in H_{k+2}(T^{k+2,n-k-2})\cong\Z$ and the standard embedding
$t=t^{2n-k-1}_{k+2,n-k-2}$;
the map $h_{t,n}$ is an isomorphism.

Construct an embedding $\tau:T^{k+1,n-k-1}\to\R^{2n-k-1}$ for $n-k$
even as follows.
Take a nonzero tangent vector field $v:S^{n-k-1}\to\R^{n-k}$ on $S^{n-k-1}$.
We have $v(a)\perp a$.
Define a map
$$\tau':\R^2\times\R^k\times S^{n-k-1}\to\R^{n-k}\times\R^k\times S^{n-k-1}
\quad\text{by}\quad\tau'(x,y,s,a)=(xa+yv(a),s,a).$$
Define an embedding $\tau$ to be the composition
of the restriction of $\tau'$ and the standard inclusion:
$$T^{k+1,n-k-1}\to\R^{n-k}\times\R^k\times S^{n-k-1}\subset\R^{2n-k-1}.$$

{\it Remark \md follows from the proof of Main Theorem \Hd.a} because the Hurewicz
homomorphism $\pi_{k+2}(N)\to H_{k+2}(N)$ is  epimorphic and by the following
easy result for $m=2n-k,2n-k-1$ (for $m=2n-k-1$ this is essentially the same
as the Summation Lemma \Pc.c).

{\it Let $N$ be a closed orientable $n$-manifold and $m\le 2n-k$.
The Whitney invariant
$W_m: E\phantom{}^m(N)\to H_{k+1}(N;\Z_{(n-k-1)})$
is defined as the composition of $W_{2n-k}$ and the map
$ E^m(N)\to E^{2n-k}(N)$ induced by the inclusion $\R^m\to \R^{2n-k}$.
If
$$f:N\to\R^{2n-k-1},\quad g:T^{k+1,n-k-1}\to\R^{2n-k-1},
\quad s:S^{k+1}\times D^{n-k-1}\to N$$
are embeddings, then $W_m(f\#_sg)=W_m(f)+W_m(g)[s|_{S^{k+1}\times 0}]$, where
$k+1<n/2$ and $W_m(g)$ is considered as an element of $\Z_{(n-k-1)}$.}


\head 3. Remarks \endhead

\smallskip
 {\bf A descriptions of generators and relations
of $ E^{2n-k-1}(S^{k+1}\times S^{n-k-1})$.}


There are 1--1 correspondences
\footnote{See [Sk08, Theorem 3.9 and table before Theorem 3.10], where
$KT^m_{p,q}= E^m(S^p\times S^q)$.
This result holds also for $n=2k+5$ in the PL category.
There is an isomorphism not only a 1--1 correspondence, for the `parametric
connected sum' group structure on the set of embeddings defined in [Sk06, \S2,
Sk08, \S3]; cf. [Sk08, \S3.4].}
$$E^{2n-1}(S^1\times S^{n-1})\to\cases
\Z_2=\left<\gamma\ |\ 2\gamma=0\right> &n=2s+1\\
\Z\oplus\Z_2=\left<\gamma,\tau\ |\ 2\gamma=0\right> &n=2s \endcases
\quad\text{for }n\ge6,\quad\text{and}$$
$$E^{2n-k-1}(S^{k+1}\times S^{n-k-1})\to\cases
0 &n-k=4s+1\\ \Z_2=\left<\gamma\ |\ 2\gamma=0\right> &n-k=4s+3\\
\Z_4=\left<\tau\ |\ 4\tau=0\right> &n-k=4s\\
\Z_2\oplus\Z_2=\left<\gamma,\tau\ |\ 2\gamma=2\tau=0\right>&n-k=4s+2
\endcases$$
for $n\ge2k+6$ and $k\ge1$.
The generators are $\gamma$ and, for $n-k$ even, $\tau$.
They are defined at the end of \S2.
The relations are $2\gamma=0$ always, and when $k\ge1$,
$\gamma=2\tau$ for $n-k=4s$, $2\tau=0$
for $n-k=4s+2$ and $\gamma=0$ for $n-k=4s+1$.

These facts follow by [Sk08, Theorem 3.9, Pa56] because the following diagram
is commutative (the map $r$ is induced by restriction, the maps
$\tau_k,\mu'',\nu''$ are defined in [Sk06]):
$$\CD \Z_2 @>>b>  E^{2n-k-1}(T^{k+1,n-k-1}) @>>r>
 E^{2n-k-1}(D^{k+1}\times S^{n-k-1})\\
@VV\cong V   @VV\tau_{k+1} V  @VV\tau_k V \\
\pi_{n-k-1}(S^{n-k-2}) @>> \mu'' > \pi_{n-k-1}(V_{n,k+2}) @>> \nu'' >
\pi_{n-k-1}(V_{n,k+1}) \endCD.$$


\smallskip
{\bf Remarks to Main Theorem \Hd.}


\smallskip
(a) {\it In the PL category} the dimension restriction of Main Theorem \Hd.a
can be relaxed to $n\ge k+5$ for $n-k$ even.
(Because this can be done in all lemmas of \S2.
For $n-k$ odd we need $n\ge2k+6$ because in the proof of the Twisting Lemma \tw
we use the Summation Lemma \Pc.)
It would be interesting to obtain an analogue of Main Theorem \Hd.a for
the PL case and $k+4\ge n\ge2k+2\ge4$ (for $n\le2k+1$ the manifold $N$ is
a homotopy sphere, cf. [Sk08, remark after Theorem 2.8]).
If $k+4\ge n\ge2k+2\ge4$, then

$\bullet$ for $n=k+4$ we have $(n,k)=(6,2)$ (then $N$ is a connected sum of
$S^3\times S^3$ by [Sm62, Theorem B]) or $(n,k)=(5,1)$ (cf. [CM]),

$\bullet$ for $n=k+3$ we have $(n,k)=(4,1)$, cf. [CS10], and

$\bullet$ $n\le k+2$ is impossible.

\smallskip
{\bf Conjecture.} If $(n,k)=(5,1)$ or $(6,2)$, then there are exact sequences
of sets
$$H_{k+1}(N;\Z_6)\overset b\to\to E^{2n-k-1}_{PL}(N)\overset{W\times W'}\to
\to H_{k+2}(N)\times H_{k+1}(N;\Z_2)\to0,$$
$$H_{k+1}(N)\times C_n^3\overset{b\times\#}
\to\to E^{2n-k-1}_{DIFF}(N)\overset{W\times W'}\to
\to H_{k+2}(N)\times H_{k+1}(N;\Z_2)\to C_{n-1}^3,$$
where $C_6^3=0$, $C_5^3\cong\Z_2$ and $C_4^3\cong\Z_{12}$. Cf. [Sk06].

\smallskip
(b) {\it On the Twisting Lemma \Tw.} In Main Theorem \hd for
$H_{k+1}(N)=0$ or $H_{k+2}(N)=0$ the $(k+2)$-parallelizability condition
can be dropped (because for $H_{k+1}(N)=0$ Main Theorem \hd is a
particular case of [Sk08, Theorem 2.13] or because for $H_{k+2}(N)=0$ we do
not need the Twisting Lemma \tw in the proof, respectively).

We conjecture that in the Twisting Lemma \tw the assumptions can be relaxed
using explicit construction [HH63, p. 133, Vr77, Lemma 6.1 and Proposition 7.1].
In particular, that the assertion of the Twisting Lemma \tw holds for $k=0$ and
$N=M\times S^{n-2}$, where $M$ is a sphere with handles.
The problem is that the construction is now performed not on the top cells,
and thus should be extended to top cells.

\smallskip
(c) {\it Parametric connected sum with $\tau$.}
Let $N$ be a closed $k$-connected $(k+1)$-parallelizable $n$-manifold
and $f_0:N\to\R^{2n-k-1}$ an embedding.
Assume that $n-k$ is even, $n\ge2k+6$ and $k\ge1$.
Let us define an action
$$b_1:H_{k+1}(N;G_{n-k})\to E^{2n-k-1}(N),\quad\text{where}
\quad G_{n-k}=\cases\Z_2 &n-k\equiv2\mod4\\  \Z_4 &n-k\equiv0\mod4\endcases.$$
Set $b_1(x):=f_0\#_{\overline x}\tau$, where
$\overline x:S^{k+1}\times D^{n-k-1}\to N$ is a smooth embedding representing
$x\in H_1(N)$.
One can prove that this is well-defined.

If $n-k=4s+2$, then $Wb_1(x)f=\rho_2(xW(\tau))=\rho_2x$, so $Wb_1=\id$.
Hence the right $H_1(N;\Z_2)$ of the exact sequence from the
Main Theorem \Hd.a is a factor of $E^{2n-k-1}(N)$.

If $n-k=4s$, then $W'b_1=\rho_2$,
$b=2b_1$ (because $b_t=2\tau$), so the following sequence of sets is exact
$$H_{k+1}(N;\Z_4)\overset{b_1}\to\to E^{2n-k-1}(N)\overset W\to\to H_{k+2}(N)\to0.$$

(d) Note that $H_{k+2}(N)$ is a factor of $ E^{2n-k-1}(N)$ if
it is a direct summand of $\pi_n(C)$ in the Complement Lemma \Co.

\smallskip
(e) {\it The kernel of the map $b'\circ AD^{-1}$ from the Complement Lemma} \co
is the image of the map $\beta:H_{n+1}(C)\to H_{n-1}(C;\Z_2)$ from the
Whitehead sequence [Wh50].
D. Crowley conjectured that $\beta$ is the composition of $\rho_2$ and the
linear dual of $\Sq^2$.

By Alexander duality, $H_{n+1}(C)\cong H_{k+3}(N)$.
So $b'\ne0$ for $N=S^{k+1}\times S^{n-k-1}$, although $b=0$ for $n-k=4s+1$.
This is so because $\psi$ need not be injective, i.e. embeddings $b(x)$ and
$b(y)$ could be isotopic although not isotopic relative to $N_0$.

The above equality $b=2b_1$ implies that $\ker b$ contains
the image of the Bockstein homomorphism $H_{k+2}(N;\Z_4)\to H_{k+1}(N;\Z_2)$
(this image could be anyway trivial because of $(k+1)$-parallelizability).

We conjecture that if $n-k\ne4s+1$, then $b'$ is injective for
$k\ge2$ and $\ker b'=\rho_2 H_4(N)\cap w_2(N)$ for $k=1$.


\smallskip
(f) {\it An idea how to reduce, for $k\ge2$, Main Theorem \Hd.a
to the cited result [BG71, Corollary 1.3].}
Denote $u=n-k-1$.
Since $N$ is $k$-connected, $N_0$ is homotopy equivalent to an $u$-polyhedron.
The space $V_{M,M-u+1}$ is $(u-2)$-connected, $\pi_{u-1}(V_{M,M-u+1})$ is $\Z$
or $\Z_2$ according to $u$ even or odd, and $\pi_{u-1}(V_{M,M-u+1})$ is
\ $\Z_4$, \ $0$, \ $\Z_2\oplus\Z_2$, \ $\Z_2$ \ according to
$u\equiv3,0,1,2\mod4$ [Pa56].
Now use the following form of homotopy classification of maps from
$u$-polyhedron to $(u-2)$-connected space, cf. [Po50, SU51].

\smallskip
{\bf Conjecture.}
{\it Let $X$ be an $u$-dimensional complex and $Y$ an $(u-2)$-connected complex.
For $u\ge4$ there are maps $W,W'$ and an action $b$ for which the following
sequence of sets is exact:
$$H^u(X;\pi_u(Y))\overset b\to\to[X;Y]\overset{W\times W'}\to
\to H^{u-1}(X;G)\times H^u(X;G\otimes \Z_2),\quad G:=\pi_{u-1}(Y).$$
For $u$ odd $W\times W'$ is surjective, and for $u$ even
$im(W\times W')=\{(p,q)\ :\ q=(\beta\otimes id_G) p\}.$}

\smallskip
(g) The following assertion allows to make Definition \bn even more explicit in the
PL category:

\smallskip
{\bf Lemma.}
{\it For each $n\ge3$ there is a PL embedding
$\eta'=\eta'_n:S^n\to S^{n-1}\times D^2$ whose composition with the projection
onto $S^{n-1}$ is homotopic to $\Sigma^{n-3}\eta$.}

\smallskip
{\it Proof.}
Indeed, define an embedding $\eta'_3:S^3\to S^2\times D^2$ by
$\eta'_3(z_1,z_2):=((z_1:z_2),z_1)$.
For $n\ge4$ define $\eta'_n$ to be the composition of $\Sigma\eta'_{n-1}$ and
the inclusion $\Sigma(S^{n-1}\times D^2)\subset S^n\times D^2$.
By induction we see that $\eta'_n$ is as required.

\bigskip
{\bf The Haefliger-Wu invariant for embeddings of $n$-manifolds into
$\R^{2n-1}$.}

We shall use the following notation for $X=N$ or $X=N_0$, although it makes
sense for general $X$.
Denote $\t X:=X\times X-\diag$.
Let $\pi_{eq}^m(\t X)$ be the set of equivariant maps
$\t X\to S^m$ up to equivariant homotopy.
For the definition of the {\it Haefliger-Wu invariant}
$\alpha: E^m(X)\to\pi^{m-1}_{eq}(\t X)$ see [Sk08, \S5].

Denote $X^*:=(X\times X-\diag)/(x,y)\sim(y,x)$.
Consider the groups
$H^i(X^*;\Z_{tw})$ with coefficients twisted according to the double cover
$q:X\times X-\diag\to X^*$ (see the details of this definition in [Ba75] where
$\Z_{tw}$ was denoted by $Zg$).

\smallskip
{\bf The Bausum Theorem.} [Ba75, Proposition 5], cf. [Ya84, To]
{\it For each $n$ and closed $n$-manifold $N$
if the set $\pi^{2n-2}_{eq}(\t N)$ is non-empty, then this set posesses an
abelian group structure such that the following sequence is exact:
$$H^{2n-3}(N^*;\Z_{tw})\overset A\to\to H_1(N;\Z_2)\overset{\t b}\to\to
\pi^{2n-2}_{eq}(\t N)\overset\deg\to\to H^{2n-2}(N^*,\Z_{tw})\to0,\quad
\text{where}$$
$$uA(\gamma)=\cases \Sq^2\gamma & n\text{ odd}\\
\Sq^2\gamma+w_1(q)^2\gamma & n\text{ even} \endcases$$
for certain isomorphism $u:H_1(N;\Z_2)\to H^{2n-1}(N^*;\Z_2)$. }

\smallskip
The assumption $n\ge6$ in [Ba75, Proposition 5] was used to apply the Haefliger
Embedding Theorem but not for the proof of the above result.
We conjecture that $\t b=\alpha b$, where $b$ is defined in \S1.

\smallskip
{\bf Deleted Product Lemma. }
{\it For $n\ge4$ and a connected orientable $n$-manifold $N_0$
with non-empty boundary there is a 1--1 correspondence
$ E^{2n-1}(N_0)\overset{\deg\circ\alpha}\to\to H^{2n-2}(N_0^*;\Z_{tw})$.}

\smallskip
{\it Proof.} The lemma follows because there are 1--1 correspondences
$$ E\phantom{}^{2n-1}(N_0)\overset\alpha\to\to\pi_{eq}^{2n-2}(\t N_0)
\overset\deg\to\to H^{2n-2}(N_0^*;\Z_{tw}).$$
The map $\alpha: E^m(N_0)\to\pi_{eq}^{m-1}(\t N_0)$ is bijective for
$2m\ge3n+2$ by [Ha63, 6.4, Sk02, Theorems 1.1.$\alpha\partial$ and
1.3.$\alpha\partial$] because $2\cdot(2n-1)\ge3n+2$ for $n\ge4$.

The map $\deg$ is 1--1 correspondence by an equivariant analogue of the
Steenrod homotopy classification theorem (which states that
{\it $\deg^{-1}(u)$ is in 1--1
correspondence with certain quotient of $H^{2n-1}(N_0^*;\Z_2)$}, cf.
[Ba75, Proposition 5, GS06, beginning of proof of the Theorem]) because
$H^{2n-1}(N_0^*;\Z_2)=0$.

(The latter isomorphism for $n=1$ this is obvious.
For $n\ge2$ we have $2\cdot2n\ge3n+2$ and $N$ is connected, hence
there are 1--1 correspondences
$H^{2n-1}(N_0^*)\to\pi_{eq}^{2n-1}(\t N_0)\overset{\alpha_{2n}^{-1}}\to\to
 E^{2n}(N_0)\to\{0\}$.
Here the first 1--1 correspondence exists by the equivariant Hopf Theorem
(which follows e.g. from [GS06, beginning of proof of the Theorem]), the second
one exists by [Ha63, 6.4, Sk02, Theorems 1.1.$\alpha\partial$ and
1.3.$\alpha\partial$] because $2\cdot(2n-1)\ge3n+2$ for $n\ge2$, and the third
one exists by [HH63, 3.1].
There is of cource a purely algebraic proof of this fact.
Then $H^{2n-1}(N_0^*;\Z_2)=
(H^{2n-1}(N_0^*)\otimes\Z_2)\oplus(H^{2n}(N_0^*)*\Z_2)=0$.)
\qed



\enddocument
\end

\ref \key Po64 \by M. M. Postnikov
\paper Investigations of the homotopy theory of continuous mappings
 \yr 1964
\jour AMS Translations, Four Papers on Topology
\endref

Classification of mappings of an $(n+2)$-complex into an $(n-1)$-connected
space with vanishing $(n+1)$-st homotopy group
Nobuo Shimada and Hiroshi Uehara
Source: Nagoya Math. J. Volume 4 (1952), 43-50.

\newpage
\head 4. The Whitney invariant for embeddings of 2-manifolds in the 4-space \endhead


Complete classification of embeddings of a given $n$-manifold $N$ into
$S^{n+2}$ up to isotopy (or concordance) seems to be hopeless because it is
such for $N=S^n$.
So it is interesting to obtain complete classification of embeddings of a
given $n$-manifold $N$ into $S^{n+2}$ `modulo knots $S^n\to S^{n+2}$'.
The notion of almost concordance (deined below) is not only useful to study the
initial
problem (of classification of embeddings up to concordance) for $m\ge n+3$, but
is a good notion of 'concordance modulo knots $S^n\to S^{n+2}$', because any
knot $S^n\to S^{n+2}$ is almost concordant to the trivial knot.
Cf. [MR].

Let $N$ be a closed connected $n$-manifold.
Two embeddings $f,g:N\to\R^m$ are called {\it almost isotopic}, if
there exists an embedding $F:N\times I\to\R^m\times I$ such that

$F(x,0)=(f(x),0)$ for each $x\in N$,

$F(x,1)=(g(x),0)$ for each $x\in N$,

$F(N\times\{t\})\subset\R^m\times\{t\}$ for each $t\in I$.

$\Sigma(F):=\{x\in N\times I\ |\ \#F^{-1}Fx>1\}\subset B^n\times I$
for some $B^n\subset N$.

In the smooth category we assume additionally that $dF$ is non-degenerate
over $(N-B^n)\times I$.

By [Sk02, Theorem 2.2.q, cf. Sk07, Corollary (a) in p. 258]

{\it for a sphere with handles (i.e. a closed connected orientable 2-manifold)
$N$ the set of PL almost embeddings $N\to\R^4$ up to PL almost concordance is
in 1--1 correspondence with $H_1(N;\Z_2)$.}

The 1-1 correspondence is the Whitney invariant $W$ defined below.

In this paper we prove that each element $H_1(N;\Z_2)$ is the Whitney invariant
of some {\it topological embedding}, not only {\it PL almost embedding}.

\smallskip
{\bf Theorem.} {\it Let $N$ be sphere with handles.
For each $\varphi\in H_1(N;\Z_2)$ there is a} topological
embedding $f:N\to\R^4$ {\it such that $W(f)=\varphi$.}

\smallskip
We also prove a version of the above result on almost embeddings.

\smallskip
{\bf Theorem 2.}
{\it Denote by $\widehat{ E}_{DIFF}^4(N)$ the set of smooth almost isotopy
classes of smooth embeddings $N\to\R^4$.
Then the Whitney invariant
$$W:\widehat{ E}\phantom{}_{DIFF}^4(N)\to H_1(N;\Z_2).$$
is injective}.

\smallskip
{\it The first definition of $w(f)$ for a TOP embedding $f:T^2\to\R^4$},
cf. [Sk02, \S6, Sk08, \S6].
Let $f:T^2\to\R^4$ be a TOP embedding.
Define a map
$$\hat f_r:S^1\times S^1\times S^1\to S^3\quad\text{by}\quad
\hat f_r(x,y,z)=\frac{f(x,y)-f(-x,z)}{|f(x,y)-f(-x,z)|}.$$
This map is equivariant with respect to the involution $(x,y,z)\to(-x,z,y)$ on
$S^1\times S^1\times S^1$ and the antipodal involution on $S^3$.
The set of such equivariant maps $S^1\times S^1\times S^1\to S^3$ up to
equivariant homotopy iz $\Z_2$ [?].
Let $w_r(f)\in\Z_2$ be the class of $\hat f_r$.

Analogously construct $w_l(f)$ from the map
$$\hat f_l:S^1\times S^1\times S^1\to S^3\quad\text{defined by}\quad
\hat f_l(x,y,z)=\frac{f(x,z)-f(y,-z)}{|f(x,z)-f(y,-z)|}.$$
{\it The Whitney invariant} is $w(f)=(w_l(f),w_r(f))$.


\smallskip
{\it The second definition of $w(f)$ for a DIFF (or PL) embedding $f:N\to\R^4$
of a sphere with handles $N$}, cf. [Sk08, \S2].
Then $w(f)$ is the homology class of the sum $[\Cl\Sigma(H)]$ of 1-dimensional
simplices of the self-intersection set $\Cl\Sigma(H)$ of a general
position homotopy $H$ between $f$ and $in$.

This $w(f)$ is indeed independent on the choice of $H$ (because for a homotopy
$X$ from $H_1$ to $H_2$ we have
$\partial[\Cl\Sigma(X)]=[\Cl\Sigma(H_1)]-[\Cl\Sigma(H_2)]$).

\smallskip
{\it The third definition of $w(f)$ for a smooth embedding $f:N\to\R^4$ of a
sphere with handles $N$.}
Let $B^2$ be a closed 2-ball in $N$. Denote $N_0:=\Cl(N-B^2)$. Let
$f:N\to S^4$ be an arbitrary smooth embedding.
Then $f|_{N_0}$ is smoothly (and hence ambiently) isotopic to $in|_{N_0}$
[HH63, 3.1.b].
Make an isotopy of $f$ and assume that $f=in$ on $N_0$.

Take a general position homotopy $F:B^2\times I\to S^4$ relative to
$\partial B^2$ between the restrictions of $f$ and $in$ to $B^2$.
Then $f\cap F:=(f|_{N-B^2})^{-1}F(B^2\times I)$ (i.e. 'the
intersection of this homotopy with $f(N-B^2)$') is a 1-manifold
(possibly non-compact) without boundary.
Define $w(f)$ to be the homology class of the closure of this oriented
1-manifold:
$$w(f):=[\Cl(f\cap F)]\in H_1(N_0,\partial N_0;\Z_2)\cong H_1(N;\Z_2).$$
Or, equivalently, $w(f):=[f|_{B^2}\cup in|_{B^2}]\in H_2(\R^4-in(N-B^2))
\cong H^1(N-B^2)\cong H_1(N-\Int B^2,\partial B^2)\cong H_1(N,B^2)
\cong H_1(N)$ (with $\Z_2$-coefficients).

It is clear that the third definition is equivalent to the first one and to
the second one.
Hence the 'third Whitney invariant' is well-defined, i.e. independent on the
choice of $F$ and the isotopy making $f=in$ outside $B^2$.

\smallskip
{\it Proof of Theorem (a).}
Cf. [Sk08', PL Classification Theorem, Sk02, Theorem 2.2.q, Sk07, Almost
Embedding Theorem (a)].
Let $f,g:N\to\R^4$ be smooth embeddings such that $w(f)=w(g)$.
We use the notation from the first paragraph of the third
definition of the Whitney invariant.
Analogously to the argument there by a smooth isotopy of $f$ we may assume that
$f=g$ on $N_0$.

Let $W$ be the closure of the complement in $S^4$ to the space of the normal
bundle of $g$ restricted to $N_0$ (i.e. to the regular neighborhood of $gN_0$
modulo $g\partial B^2$).
By general position $W$ is simply connected.
Hence the Hurewicz homomorphism $\pi_2(W)\to H_2(W;\Z)$ is an isomorphism.
By Alexander and Poincar\'e duality
$H_2(W)\cong H^1(N_0)\cong H^1(N)\cong H_1(N)$.
The class $[f|_{B^2}\cup g|_{B^2}]\mod2$ goes to $w(f)-w(g)$ under this
isomorphism.

By a smooth isotopy of $f$ (not fixed on $N_0$) we can change
$[f|_{B^2}\cup g|_{B^2}]\in\pi_2(W)$ by any even element (the homotopy
$f_t:T^2\to\R^3$ which is projection of such an isotopy $f_t:T^2\to\R^4$ is
shown in Fig. 1).
So if $w(f)=w(g)$, we may assume that $[f|_{B^2}\cup g|_{B^2}]=0\in\pi_2(W)$.
Hence $f|_{B^2},g|_{B^2}:B^2\to W$ are homotopic fixed on $\partial B^2$.
Therefore $f$ and $g$ are smoothly almost isotopic.
\qed

\smallskip
{\it Proof of Theorem (b).}
We prove the result for $N=T^2$, for the general case the proof is analogous.

Consider the standard embedding $2D^3\times S^1\subset\R^4$.
The {\it right Hudson torus} $\tau_r$ is the connected sum of
$2\partial D^3\times x$ with
$\partial D^2\times S^1\subset D^3\times S^1\subset 2D^3\times S^1\subset\R^4$.
The {\it left Hudson torus} $\tau_l$ is defined analogously but interchanging
factors.
Denote by $in:T^2\to\R^4$ the standard embedding.
One can check that $w(in)=0$, $w(\tau_r)=(0,1)$ and $w(\tau_l)=(1,0)$.
This is most clear by the third definition of $w(f)$ below.

Denote by $W$ the closure in $S^4$ of the complement of $S^4$ to a regular
neighborhood in $S^4$ of $in(T^2-\Int B^2)$ modulo $in(\partial B^2)$.
Then $W\simeq S^4-(S^1\vee S^1)\simeq S^2\vee S^2$.
Hence $H_2(W,\partial W;\Z)\cong H^2(W;\Z)\cong\Z\oplus\Z$.
Identify these groups.
Since $W$ is simply-connected, the embedding $in|_{\partial B^2}$ can be
extended to a map $\varphi:B^2\to W$.
We may assume that $\varphi$ is proper.
Since the Hurewicz homomorphism $\pi_2(W)\to H_2(W;\Z)$ is an isomorphism, we
may assume that $\varphi$ represents $(1,1)\in\Z\oplus\Z=H_2(W,\partial W;\Z)$.
Then $\varphi$ is homotopic to an embedding $\overline\varphi$ by
[GM86, Casson's lecture I, Embedding Theorem and remark in brackets] (because
$(1,1)\cap(1,0)=1$).
Now define $f:T^2\to\R^4$ by $f(x)=x$ for $x\in T^2-B^2$ and
$f(x)=\overline\varphi(x)$.

Analogously to the proof of equivalence of the third and the first
definitions of the Whitney invariant we obtain that $w(f)=(1,1)$.
\qed

\smallskip
{\bf Conjecture.} {\it The set $E^4(T^2)$ of smooth (or PL, or TOP locally flat)
embeddings
$T^2\to\R^4$ up to smooth almost isotopy consists of four elements.

Or, equivalently, there is a smooth (or PL, or TOP locally flat) embedding
$T^2\to\R^4$ that is not
smoothly (or PL, or TOP locally flat) almost
isotopic either to the standard embedding, or to the right Hudson torus,
or to the left Hudson torus (see below).


Or, equivalently, there is a smooth (or PL, or TOP locally flat) embedding
as in Theorem (b). (In particular, there is a smooth embedding
$f:T^2\to\R^4$ such that $w(f)=(1,1)$.)}


\smallskip
Analogues of Conjecture are apparently true in the PL category and/or for
sphere with handles $N$.
Analogous results can be stated for closed non-orientable manifolds.
If smooth embeddings $N\to\R^4$ are PL almost isotopic, then they are DIFF
almost isotopic (because a concordance can be smoothed outside a point).

\bigskip
{\bf Some further remarks and conjectures.}

\smallskip
Note that there is a bijection $H_1(N;\Z_2)\to\pi^3_{eq}(\t N)$ mapping the
Whitney invariant to the Haefliger-Wu $\alpha$-invariant.

Corollary: {\it Almost concordant (even topologically) smooth embeddings
of sphere with handles into $\R^4$ are smoothly almost isotopic}.

S. Melikhov: To, chto gladkaja pochti-izotopija ekvivalentna topologicheskoj
ochevidno: na (odnomernom) spine'e sfery s ruchkami ona gladkaja po obshemu
polozheniju, na ego okrestnosti mozhet byt' sdelana gladkoj, a dopolnenie k
etoj okrestnosti - disk. Analogichno dlja pochti-konkordantnosti.

Samo sledstvie ne imeet ko mne nikakogo otnoshenija i trivialno
dokazyvaetsja standartnoj tehnikoj: konkordantnost' raskladyvaetsja na
ruchki (nachinaja s ejo ogranichenija na spine), dalee ruchki
raspolagajutsja po vysote v porjadke vozrastanija indeksov i
sokrashajutsja. Snachala vyprjamljaetsja ogranichenie konkordantnosti na
spine'e, dalee ruchki indexa 2 peredvigajutsja vverh tak chtoby ne peresech'
konkordantnost' spine'a, a mezhdu soboj im mozhno peresekat'sja.

Corollary follows from Theorem (a) below because $w(f)$ is invariant under
topological almost concordance of $f$ (by the first definition of $w(f)$ below
and the Suspension Theorem).

\smallskip
(2) {\it Action of the group of self-homeomorphisms of the torus on
on $H_1(T^2;\Z_2)$.}

If $A$ is the linear autohomeomorphism of $T^2$ that exchanges parallel
and meridian, then $in\circ A$ is trivial.
The same if $A$ is the Dehn twist by 2 along the meridian of $T^2$.
If $A$ is the Dehn twist by 1 along the meridian of $T^2$, then $in\circ A$
is not isopositioned to the standard embedding [Mo83].

Clearly, $w(\tau\circ A)=w(in\circ A)+A^{-1}_*w(\tau)$ for any embedding
$\tau:T^2\to\R^4$ and any automorphism $A:T^2\to T^2$.
If (a) in Lemma below holds, then $w(\tau_r\circ X)=(1,1)$ and the
Conjecture is wrong.
Note that $SL_2(\Z_2)\cong S_3$.

\smallskip
{\bf Lemma.} {\it Denote $U=\left(\matrix 0 & 1 \\ 1 & 0 \endmatrix\right)$
and $X=\left(\matrix 1 & 1 \\ 1 & 0 \endmatrix\right)$.
Suppose that $w:SL_2(\Z_2)\to\Z_2\oplus\Z_2$ is a map such that $w(E)=w(U)=0$
and $w(AB)=w(B)+B^{-1}w(A)$ for each $A,B\in SL_2(\Z_2)$.
Then either

(a) $w(A)=0$ for each $A\in SL_2(\Z_2)$, or

(b) $w(X)=w(UX)=(0,1)$ and $w(XU)=w(UXU)=(1,0)$. }

\smallskip
{\it Proof.} Since $w(U)=0$, it follows that $w(AU)=Uw(A)$ and $w(UA)=w(A)$.
Then $w(UAU)=Uw(A)$.
Recall that $UA$, $AU$ and $UAU$ is obtained from $A$ by transposition of
lines, by transposition of columns by central symmetry, respectively.

We have $(UX)^2=E$, so $X^{-1}=UXU$, hence $X^2=UXU$.
Therefore $Uw(X)=w(X^2)=w(X)+X^{-1}w(X)$, so $(U+E+UXU)w(X)=0$, thus
$w(X)=(0,\alpha)$.
Then $w(UX)=(0,\alpha)$ and $w(XU)=w(UXU)=(\alpha,0)$.
\qed


\smallskip
(3) Another equivalence relation 'cancelling' knots is as follows.
Two embedding $f_1,f_2:T^2\to\R^4$ are called {\it PL isotopic}, if there is a
knot $g:S^2\to\R^4$ such that $f_1$ is ambient isotopic to $f_2\#g$
($f_2T^2$ and $gS^2$ are contained in disjoint 3-balls).
We conjecture that {\it there are non-equivalent but almost concordant
embeddings $T^2\to\R^4$}.
In order to prove the conjecture one could try to construct a non-trivial
embedding $f:T^2\sqcup S^2\to\R^4$ (possibly using $\beta$-invariant), and then
try to prove that ('linked') connected sum of components of $f$ is not
equivalent to $f|_{T^2}$.

It is clear that for any embedding $g:S^2\to\R^4$ we have $w(f\# g)=w(f)$.
 Is the Hudson torus isotopic to a connected sum of the
standard torus with a knot?

(4) For $u\in\Z$ instead of an embedded 2-sphere $2\partial
D^3\times x$ we can take $u$ copies $(1+\frac1n)\partial D^3\times
x$ ($n=1,\dots,u$) of 2-sphere outside $D^3\times S^1$ 'parallel'
to $\partial D^3\times x$. Then we join these spheres by tubes so
that the homotopy class of the resulting embedding $S^2\to
S^4-D^3\times S^1\simeq S^4-S^1\simeq S^2$ will be $u$. Let $f_u$
be the connected sum of this embedding with the standard embedding
$\partial D^2\times S^1\subset S^4$.
Then $w(f_u)=u\mod2$.

Is $f_2$ ambient isotopic to $f_0$?

(5) Another invariant of an embedding $f:T^2\to\R^4$ is {\it the Rokhlin
quadratic form $q(f):H_1(T^2;\Z_2)\to\Z_2$} [Mo83, Hi].
Note that there are exactly 4 bilinear forms
$H_1(T^2;\Z_2)\times H_1(T^2;\Z_2)\to\Z_2$ and so probably 4 quaddratic
forms $H_1(T^2;\Z_2)\to\Z_2$.
It would be interesting to find a relation of $q(f)$ and $w(f)$.

(6) An analogue of the Whitney invariant can be defined for embeddings
$S^1\sqcup S^1\to\R^3$. It would be equal to the linking coefficient.

\enddocument